\documentclass[11pt]{article}
\usepackage[utf8]{inputenc}
\usepackage[T1]{fontenc}
\usepackage{comment}
\usepackage{authblk}
\usepackage{relsize}
\usepackage{changepage}
\usepackage{mathtools}
\usepackage[english]{babel}
\usepackage{amsmath}
\usepackage{amsthm}
\usepackage{amsfonts}
\usepackage{amssymb}
\usepackage{graphicx}
\usepackage[usenames,dvipsnames]{color}
\theoremstyle{plain}
\newtheorem{theorem}{Theorem}[section]
\newtheorem{corollary}[theorem]{Corollary}

\newtheorem{lemma}[theorem]{Lemma}

\newtheorem{observation}[theorem]{Observation}

\makeatletter
\newcommand{\vast}{\bBigg@{4}}
\newcommand{\Vast}{\bBigg@{5}}
\makeatother
\definecolor{bulgarianrose}{rgb}{0.28, 0.02, 0.03}
\definecolor{gray}{rgb}{0.5, 0.5, 0.5}

\theoremstyle{definition}

\newtheorem{conjecture}[theorem]{Conjecture}

\theoremstyle{remark}

\usepackage[left=2cm,right=2cm,top=2cm,bottom=2cm]{geometry}
\usepackage{amssymb}
\makeatletter
\def\namedlabel#1#2{\begingroup
    #2%
    \def\@currentlabel{#2}%
    \phantomsection\label{#1}\endgroup
}
\usepackage[normalem]{ulem}
\usepackage{dsfont}
\usepackage{accents}
\usepackage{verbatim}
\usepackage{ulem}
\usepackage{pgf,tikz,pgfplots}
\pgfplotsset{compat = 1.16}
\usepackage[all]{xy} 

\usepackage{stackengine}
\stackMath
\newcommand\tsup[2][2]{%
 \def\useanchorwidth{T}%
  \ifnum#1>1%
    \stackon[-.5pt]{\tsup[\numexpr#1-1\relax]{#2}}{\scriptscriptstyle\sim}%
  \else%
    \stackon[.5pt]{#2}{\scriptscriptstyle\sim}%
  \fi%
}

\usepackage{subfig}

\usepackage{enumerate}

\title{\scshape
  On the chromatic number of the preferential attachment graph}

\author{Lyuben Lichev}
\affil{Ecole Normale Supérieure de Lyon, Lyon, France}

\begin{document}

\maketitle
 
\begin{abstract}
We prove that for every $m\in \mathbb N$ and every $\delta\in (-m,0)$, the chromatic number of the preferential attachment graph $PA_t(m, \delta)$ is asymptotically almost surely equal to $m+1$. The proof relies on a combinatorial construction of a family of digraphs of chromatic number $m+1$ followed by a proof that asymptotically almost surely there is a digraph in this family, which is realised as a subgraph of the preferential attachment graph.
\end{abstract}

\hspace{1em}Keywords: preferential attachment, chromatic number

\hspace{1em}MSC Class: 05C80, 05C15

\section{Introduction}

The chromatic number is one of the most basic characteristics of any graph. For a graph $G$, the \textit{chromatic number} of $G$, denoted by $\chi(G)$, is the minimal number of colours needed to colour the vertices of $G$ so that no monochromatic edges are present. Any colouring of this type is called \textit{proper}. Some of the first instances of the graph colouring problem were in connection with colouring planar graphs in the form of maps. This led to the celebrated Four Colours Theorem stating that any planar graph may be properly coloured in four colours. The literature on graph colouring problems and their applications has by now become nearly spanless.\par

Another line of research was initiated by Erd\H{o}s and Rényi in their seminal papers \cite{ER_1959} and \cite{ER_1960}, which introduced the models $G(n,p)$ and $G(n,M)$ for the first time. Later Bollob\'as proved in \cite{Bol_1988} that asymptotically almost surely
\begin{equation*}
\dfrac{n\log\left(\frac{1}{1-p}\right)}{2\log(n)}\left(1+\dfrac{\log(\log(n))}{\log(n)}\right)\le \chi(G(n,p))\le \dfrac{n\log\left(\frac{1}{1-p}\right)}{2\log(n)}\left(1+\dfrac{3\log(\log(n))}{\log(n)}\right).
\end{equation*}
Moreover, in the regime $p = n^{-\alpha}$ for every $\alpha\in (1/2, 1)$ Bollob\'as proved that $\chi(G(n,p))$ takes one of only five values asymptotically almost surely. See also \cite{Bol_book}.\par

Another random graph model - the configuration model - was introduced in~\cite{BC} and further developed by Bollobás in~\cite{Bol_book} and by Wormald in~\cite{W78}. This model is in tight connection with the uniform distribution on random graphs with bounded degree sequence. Indeed, in the case of random $d$-regular graphs with $d$ fixed for example, the uniform distribution on this set of graphs is contiguous with the $d$-regular configuration model, see \cite{Wor_1999}. Concerning the chromatic number of random regular graphs, Achlioptas and Moore show in \cite{AM_2004} that for every $d\ge 3$, the chromatic number of the uniform random $d$-regular graph is asymptotically almost surely among $k,k+1$ and $k+2$, where $k$ is the smallest integer such that $d < 2k\log(k)$. Later Coja-Oghlan, Efthymiou and Hetterich prove in \cite{C-OEH_2015} that for some $k_0\in \mathbb N$ and for every $k > k_0$ there exists a value $d_{k, col}\in \mathbb N$, for which, if $d > d_{k, col}$, the random $d$-regular graph cannot be coloured in $k$ colours  asymptotically almost surely and if $d < d_{k, col}$, such a colouring exists asymptotically almost surely.\par
 
A more recent though by now classical random graph model is the preferential attachment model (or the PA model for short) introduced by Albert and Barab\'asi in \cite{AB_1999} and \cite{AB_2002} and further formalised by Bollob\'as and Riordan in \cite{BR_2003}. This model was thouroughly studied since it shares many important characteristics with real world networks. Many technological, biological and social networks have been empirically proved to behave like the preferential attachment model with given parameters.\par

There is a number of variants of the PA model. In this paper, we follow the definition given in \cite{GS}. Fix $m\ge 1$ and $\delta > -m$. Then $(PA_t(m, \delta))_{t\ge 1}$ is a sequence of random graphs defined as follows:
\begin{itemize}
    \item for $t = 1$, $PA_1(m, \delta)$ consists of a single vertex with no edges;
    \item for $t = 2$, $PA_2(m, \delta)$ consists of two vertices with $m$ edges between them;
    \item for every $t\ge 3$, $PA_t(m, \delta)$ is defined from $PA_{t-1}(m, \delta)$ by adding one vertex $v_t$ and $m$ edges $e^1_t, e^2_t, \dots, e^m_t$ in a consecutive fashion so that, for every $j\in [m]$ and $i\in [t-1]$, the probability that, $e^j_t = v_tv_i$ is given by
    \begin{equation*}
        \dfrac{d_{t-1, j-1}(v_i)+\delta}{2m(t-2)+(t-1)\delta+j-1},
    \end{equation*}
    where $d_{t-1, j-1}(v_i)$ is the degree of the vertex $v_i$ in the graph $PA_{t-1}(m, \delta)\cup (e_t^r)_{r\in [j-1]}$.
\end{itemize}

The parameter $\tau = 3 + \dfrac{\delta}{m}$ is perhaps the most important characteristic of the model since it describes the degree distribution in $PA_t(\delta, m)$. In particular, the fraction of vertices of degree $k$ goes as $(1+o_k(1))k^{-\tau}$ asymptotically as $k\to +\infty$. One has empirically observed that in most real world networks $\tau\in (2,3)$ and thus most effort has been directed to understanding the model in this particular regime.\par

From the point of view of the graph colouring problem, the PA model was studied by Kovalenko in \cite{Kov}, who showed that for every $\varepsilon > 0$ there is $m(\varepsilon)\in \mathbb N$ such that for every $m\ge m(\varepsilon)$ asymptotically almost surely
\begin{equation*}
    \dfrac{m}{(4+\varepsilon)\log(m)}\le \chi(PA_t(m, \delta)).
\end{equation*}
He also gave a deterministic upper bound:
\begin{equation*}
    \chi(PA_t(m, \delta)) \le \dfrac{\log(t)-\log(m)}{\log\left(1+\frac{1}{m}\right)}+m +1.
\end{equation*}
His approach was based on proving lower and upper bounds on the size of the maximal independent set in the PA graph. This was carried out in a more general class of preferential attachment models, which in particular include the one described above as a special case.\par

The main result of this paper improves this result by proving that the chromatic number of the PA graph converges almost surely to $m+1$.

\begin{theorem}\label{main th}
For every $m\in \mathbb N$ and $\delta\in (-m,0)$ the chromatic number of the preferential attachment graph $(PA_t(m, \delta))_{t\ge 1}$ almost surely converges to $m+1$.
\end{theorem}

The proof is trivial in the case of $m=1$ since in this case $PA_t(1, \delta)$ is a non-empty tree for every $t\ge 1$ and consequently has chromatic number equal to 2. From this point we work with $m\ge 2$.\par

Our approach is largely influenced by the paper \cite{GS} of Garavaglia and Stegehuis. There the authors find the order of the expected value of the number of copies of a given graph $H$ in $PA_t(m, \delta)$ and prove that in some cases one has concentration around the expected value. In particular, they compute the expected number of triangles and prove concentration around this expected value as $t\to +\infty$. Rather than using their results directly, we carefully analyse the proofs of Theorem 2.2 and Lemma 6.1 from the paper to obtain some important information about a family of concrete graphs, which we construct in the sequel.\par

Our main probabilistic tool in this paper will be the simple but highly efficient second moment method. It is based on the following lemma, which is also called Paley-Zygmund's inequality.

\begin{lemma}[\cite{PZ_1932}, Lemma 19, page 192]\label{PZ ineq}
For every positive random variable $Z$ with finite second moment one has
\begin{equation*}
    \mathbb P(Z > 0)\ge \dfrac{\mathbb E[Z]^2}{\mathbb E[Z^2]}.
\end{equation*}
\end{lemma}

\textbf{Notation:} In this paper we denote by $\mathbb N$ the set of positive integers. For every $n\ge 1$, we denote by $[n]$ the set of integers $\{1,2,\dots,n\}$. Unless explicitly stated otherwise, graphs and digraphs will be denoted by $G,H$, vertices by $u,v,w$ and edges (oriented or not) by $e$, possibly with some lower and upper indices. Edges may also be denoted as concatenation of two vertices, for example $uv$, where the order of the vertices is important only if the edge is oriented. We denote by $\log$ the logarithm in base $e$. For a (di-)graph $G$, we denote by $V(G)$ the vertex set of $G$ and by $E(G)$ the edge set of $G$. For a vertex $v$ in a digraph $G$, the \textit{degree of $v$ in $G$} is the sum of the indegree and the outdegree of $v$ in $G$. Standard asymptotic notations like $o(\cdot), O(\cdot), \Theta(\cdot)$ are used throughout the paper. When it is not completely clear which is the limit variable, we add it to the notation and write for example $o_t(\cdot), O_t(\cdot), \Theta_t(\cdot)$. \par

\textbf{Organisation of the paper:} In Section~\ref{sec: construction} we construct a family of digraphs of chromatic number $m+1$ and outdegree bounded above by $m$ for every $m\ge 2$. In Section~\ref{sec: GS result} we analyse in detail the proof of Theorem 2.2 from \cite{GS} to find the correct order of the probability to find some of the digraphs, constructed in Section~\ref{sec: construction}, as a subgraph of $PA_t(m, \delta)$ for every $m\ge 2$ and $\delta\in (-m, 0)$. In Section~\ref{sec: proof} we prove Theorem~\ref{main th 1}, which together with Observation~\ref{ob 2.4} implies Theorem~\ref{main th} as a direct corollary. We finish with a discussion and open problems in Section~\ref{sec: conclusion}.

\section{The construction}\label{sec: construction}

In this paper we consider digraphs $G = (V, \vec{E})$ equipped with a bijective map $\sigma = \sigma_G: V(G)\to [|V(G)|]$ with the property that for every directed edge $uv\in \vec{E}(G)$ we have that $\sigma(u) > \sigma(v)$. We call any such map $\sigma$ an \textit{ordering} of $G$. 

Note that, given a set $V$ and a bijective map $\sigma: V\to [|V|]$, every digraph with vertex set $V$ and ordering $\sigma$ may be constructed as follows. Denote $n = |V|$ and fix a sequence of integers $(d_1, \dots, d_n)$ such that, for every $i\in [n]$, $d_i\in [0, i-1]$. Let $G_0$ be the null graph (i.e. the unique graph with 0 vertices). For every $i\in [n]$, add the vertex $\sigma^{-1}(i)$ to $G_{i-1}$ and construct $d_i$ edges going out of $\sigma^{-1}(i)$ to $d_i$ different vertices of $G_{i-1}$ arbitrarily. Then, $G = G_n$ is a digraph with ordering $\sigma$ in which the outdegree of the vertex $\sigma^{-1}(i)$, denoted $d^{(out)}_G(\sigma^{-1}(i))$ below, is $d_i$.

We call the graph obtained from a digraph $G$ by forgetting the orientation of every edge of $G$ the \textit{undirected version of $G$}. The chromatic number of a digraph $G$ is the chromatic number of its undirected version. A digraph $G$ is \textit{weakly connected} if its undirected version is a connected graph.

\begin{lemma}
The chromatic number of a digraph $(G, \sigma)$ is at most $\max_{i\in [n]} d^{(out)}_G(\sigma^{-1}(i)) + 1$.
\end{lemma}
\begin{proof}
We identify colours with positive integers. Colour the vertices of the digraph $(G, \sigma)$ greedily in increasing order with respect to $\sigma$. When constructing $G_i$, the vertex $\sigma^{-1}(i)$ is incident to $d^{(out)}_G(\sigma^{-1}(i))$ other vertices in $G_{i-1}$ and therefore there exists a colour in $[d^{(out)}_G(i) + 1]$, which may be used to colour the vertex $\sigma^{-1}(i)$ properly in $G_i$. It remains to note that for every $i\in [n]$ this procedure extends a proper colouring of $G_{i-1}$ to a proper colouring of $G_i$, which means that there is a proper colouring of $G = G_n$ with at most $\max_{i\in [n]} d^{(out)}_G(\sigma^{-1}(i)) + 1$ colours.
\end{proof}

\begin{corollary}\label{cor upper bound}
For every $t, m\in \mathbb N$ and $\delta > -m$ we have $\chi(PA_t(m, \delta))\le m+1$.\qed
\end{corollary}

For two digraphs $(G_1, \sigma_1)$ and $(G_2, \sigma_2)$ with $G_1\subseteq G_2$, we say that the ordering $\sigma_1$ is a \textit{restriction} of the ordering $\sigma_2$ (on $G_1$) and the ordering $\sigma_2$ is an \textit{extension} of the ordering $\sigma_1$ (on $G_2$) if for every vertex $v\in V(G_1)$ we have $\sigma_1(v) = \sigma_2(v)$. Also, for two digraphs $(G_1, \sigma_1)$ and $(G_2, \sigma_2)$ we say that the orderings $\sigma_1$ and $\sigma_2$ \textit{agree} (on $G_1\cap G_2$) if for every $u,v\in V(G_1\cap G_2)$ we have that $\sigma_1(u) > \sigma_1(v)$ if and only if $\sigma_2(u) > \sigma_2(v)$. In particular, if $\sigma_1$ is a restriction or an extension of $\sigma_2$, the two orderings agree, but the converse is not always true.\par

Now, for every $m\ge 2$ we construct a digraph of chromatic number $m+1$, in which the outdegree of each of its vertices is at most $m$. We remark that the constructed graph is triangle-free. Fix $m\ge 2$ and define a sequence of digraphs with orderings $(H_0, \pi_0), (H_1, \pi_1), \dots, (H_{m-2}, \pi_{m-2})$ as follows. Define 
\begin{equation*}
    H_0 = ((v_i)_{1\le i\le 7}, \{v_1v_2, v_2v_3, v_3v_4, v_4v_5, v_5v_6, v_6v_7, v_1v_7\}),\hspace{0.5em} \pi_0: j\in [7]\mapsto v_j\in \{v_i\}_{1\le i\le 7}.
\end{equation*}
In other words, $H_0$ is an oriented copy of 
\begin{equation*}
    C_7 = (\pi_0^{-1}(1),\pi_0^{-1}(2),\pi_0^{-1}(3),\pi_0^{-1}(4),\pi_0^{-1}(5),\pi_0^{-1}(6),\pi_0^{-1}(7)).
\end{equation*}
Then, for every $i\in [m-2]$, to define the digraph $(H_i, \pi_i)$, consider $i+2$ disjoint copies $(H^{j}_{i-1}, \pi^j_{i-1})_{j\in [i+2]}$ of the digraph $(H_{i-1}, \pi_{i-1})$. Then, for any $(i+2)-$tuple of vertices $(v^1, v^2, \dots, v^{i+2})\in V(H^{1}_{i-1})\times V(H^2_{i-1})\times \dots\times V(H^{i+2}_{i-1})$, add a vertex $w = w(v^1, v^2, \dots, v^{i+2})$ and directed edges $wv^1, wv^2, \dots, wv^{i+2}$. This forms the graph $H_i$. Finally, choose an arbitrary bijective map $\nu_i: V(H_i\setminus (\cup_{j\in [i+2]} H^j_{i-1}))\to [|V(H_i\setminus (\cup_{j\in [i+2]} H^j_{i-1}))|]$ and define

\begin{equation*}
\pi_i(v) = 
\begin{cases}
& \pi^j_{i-1}(v) + (j-1)|V(H_{i-1})|, \text{ if }v\in H^j_{i-1} \text{ for some }j\in [i+2],\\
& \nu_i(v) + (i+2)|V(H_{i-1})|, \text{ if }v\in H_i\setminus H_{i-1}.
\end{cases}
\end{equation*}

Then $\pi_i$ is an ordering of the vertex set of $H_i$, which induces the orientation of the edges of $H_i$ given above.

\begin{lemma}\label{chrom nb = m+1}
For every $m\ge 2$ and $i\in 0\cup [m-2]$, the chromatic number of $H_i$ is $i+3$.
\end{lemma}
\begin{proof}
Identify colours with positive inregers. We prove the statement by induction on $i \in 0\cup [m-2]$. For $i=0$ the statement is true since $\chi(C_7) = 3$. Suppose that the claim holds for $i-1$ for some $i\in [m-2]$.\par 

Fix a proper colouring of $H_i$. We show that it uses at least $i+3$ colours. If the graph $H^{1}_{i-1}\cup H^2_{i-1}\cup \dots \cup H^{i+2}_{i-1}$ is coloured in $i+3$ colours, we are done. Otherwise, by the induction hypothesis $\chi(H_{i-1}) = i+2$ and therefore $H^{1}_{i-1}\cup H^2_{i-1}\cup \dots \cup H^{i+2}_{i-1}$ is coloured in exactly $i+2$ colours, say the ones in $[i+2]$. Moreover, every subgraph of $H_i$ among $H^{1}_{i-1}, H^2_{i-1}, \dots, H^{i+2}_{i-1}$ contains a vertex in colour $j$ for every $j\in [i+2]$. We conclude that there is a vertex $w = w(v^1, v^2, \dots, v^{i+2})$ in $H_i\setminus (H^{1}_{i-1}\cup H^2_{i-1}\cup \dots \cup H^{i+2}_{i-1})$, whose neighbours $v^1, v^2, \dots, v^{i+2}$ are coloured in different colours. Thus, $\chi(H_i) \ge i+3$.\par 

On the other hand, colouring each of the graphs $H^1_{i-1}, H^2_{i-1}, \dots, H^{i+2}_{i-1}$ in the colours $1,2,\dots,i+2$ and colouring all remaining vertices in colour $i+3$ gives a proper colouring of $H_i$. Thus, $\chi(H_i) = i+3$. The induction is finished and the lemma is proved.
\end{proof}

Now, for a digraph $(H, \pi)$ we define the digraph $S_n(H, \pi) = G$ as the digraph formed from $H$ by adding, for every two (distinct) vertices $v', v''$ in $H$, $n$ vertices $v_1, v_2, \dots, v_n$ and $2n$ directed edges $v_1v', v_2v', \dots, v_nv'$, $v_1v'', v_2v'', \dots, v_nv''$. Informally this construction could be translated as adding $n$ "cherries" for any pair of vertices of $H$. Notice that $G$ may be equipped with an ordering $\sigma_{n, H}$, which extends $\pi$ on $H$, which is constructed as follows. Fix an arbitrary bijective mapping $\mu_n: V(G\setminus H) \to [n|V(H)|]$. We define
\begin{equation*}
    \sigma_{n, H}(v) = 
    \begin{cases}
    & \pi(v), \text{ if }v\in V(H),\\
    & \mu_n(v) + |V(H)|, \text{ if }v\in V(G\setminus H).
    \end{cases}
\end{equation*}

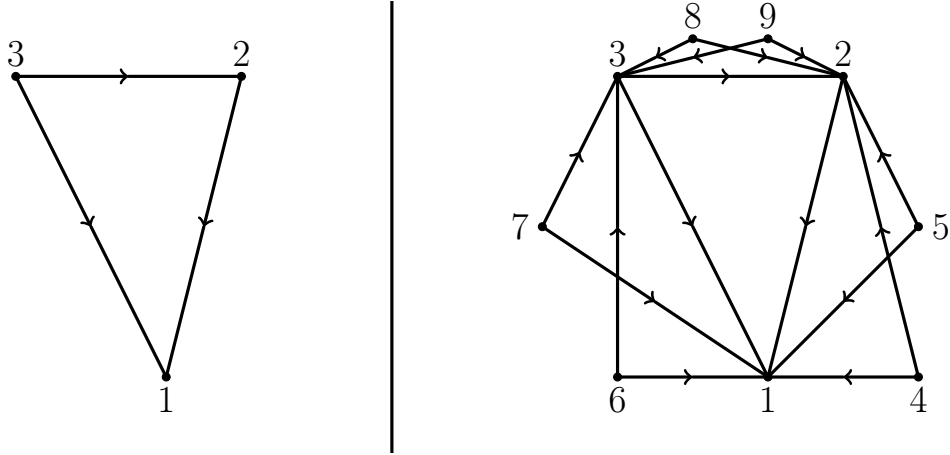
\begin{figure}
\centering
\begin{tikzpicture}[line cap=round,line join=round,x=1cm,y=1cm]
\clip(-1,-2) rectangle (13,4);

\draw [->,line width=1.2pt] (0,3) -- (1,1);
\draw [line width=1.2pt] (1,1) -- (2,-1);

\draw [->,line width=1.2pt] (0,3) -- (1.5,3);
\draw [line width=1.2pt] (1.5,3) -- (3,3);

\draw [->,line width=1.2pt] (3,3) -- (2.5,1);
\draw [line width=1.2pt] (2.5,1) -- (2,-1);

\draw [->,line width=1.2pt] (8,3) -- (9.5,3);
\draw [line width=1.2pt] (9.5,3) -- (11,3);

\draw [->,line width=1.2pt] (8,3) -- (9,1);
\draw [line width=1.2pt] (9,1) -- (10,-1);

\draw [->,line width=1.2pt] (11,3) -- (10.5,1);
\draw [line width=1.2pt] (10.5,1) -- (10,-1);

\draw [->,line width=1.2pt] (7,1) -- (7.5,2);
\draw [line width=1.2pt] (7.5,2) -- (8,3);

\draw [->,line width=1.2pt] (9,3.5) -- (8.5,3.25);
\draw [line width=1.2pt] (8.5,3.25) -- (8,3);

\draw [->,line width=1.2pt] (9,3.5) -- (10,3.25);
\draw [line width=1.2pt] (10,3.25) -- (11,3);

\draw [->,line width=1.2pt] (10,3.5) -- (10.5,3.25);
\draw [line width=1.2pt] (10.5,3.25) -- (11,3);

\draw [->,line width=1.2pt] (10,3.5) -- (9,3.25);
\draw [line width=1.2pt] (9,3.25) -- (8,3);

\draw [->,line width=1.2pt] (12,1) -- (11.5,2);
\draw [line width=1.2pt] (11.5,2) -- (11,3);

\draw [->,line width=1.2pt] (8,-1) -- (9,-1);
\draw [line width=1.2pt] (9,-1) -- (10,-1);

\draw [->,line width=1.2pt] (7,1) -- (8.5,0);
\draw [line width=1.2pt] (8.5,0) -- (10,-1);

\draw [->,line width=1.2pt] (8,-1) -- (8,1);
\draw [line width=1.2pt] (8,1) -- (8,3);

\draw [->,line width=1.2pt] (12,-1) -- (11,-1);
\draw [line width=1.2pt] (11,-1) -- (10,-1);

\draw [->,line width=1.2pt] (12,-1) -- (11.5,1);
\draw [line width=1.2pt] (11.5,1) -- (11,3);

\draw [->,line width=1.2pt] (12,1) -- (11,0);
\draw [line width=1.2pt] (11,0) -- (10,-1);

\draw [line width=1.2pt] (5,-10) -- (5,10);

\begin{scriptsize}
\draw [fill=black] (0,3) circle (1.5pt);
\draw[color=black] (0,3.3) node {\Large 3};

\draw [fill=black] (2,-1) circle (1.5pt);
\draw[color=black] (2,-1.3) node {\Large 1};

\draw [fill=black] (3,3) circle (1.5pt);
\draw[color=black] (3,3.3) node {\Large 2};

\draw [fill=black] (0+8,3) circle (1.5pt);
\draw[color=black] (0+8,3.3) node {\Large 3};

\draw [fill=black] (2+8,-1) circle (1.5pt);
\draw[color=black] (2+8,-1.3) node {\Large 1};

\draw [fill=black] (3+8,3) circle (1.5pt);
\draw[color=black] (3+8,3.3) node {\Large 2};

\draw [fill=black] (7,1) circle (1.5pt);
\draw[color=black] (6.7,1) node {\Large 7};

\draw [fill=black] (9,3.5) circle (1.5pt);
\draw[color=black] (9,3.5+0.3) node {\Large 8};

\draw [fill=black] (10,3.5) circle (1.5pt);
\draw[color=black] (10,3.5+0.3) node {\Large 9};

\draw [fill=black] (12,1) circle (1.5pt);
\draw[color=black] (12.3,1) node {\Large 5};

\draw [fill=black] (8,-1) circle (1.5pt);
\draw[color=black] (8,-1.3) node {\Large 6};

\draw [fill=black] (12,-1) circle (1.5pt);
\draw[color=black] (12,-1.3) node {\Large 4};

\end{scriptsize}
\end{tikzpicture}
\caption{On the left: a digraph $H$ with an ordering $\pi$. On the right: the digraph $G = S_2(H)$ equipped with one possible ordering $\sigma$, which coincides with $\pi$ on $V(H)$.}
\label{fig 1}
\end{figure}

\noindent See Figure~\ref{fig 1}. 

\begin{observation}\label{ob 2.4}
For every $n\ge 1$ and $m\ge 2$, the chromatic number of $S_n(H_{m-2}, \pi_{m-2})$ is $m+1$.
\end{observation}
\begin{proof}
For all $n\ge 1$ and $m\ge 2$ we have that, on the one hand, $H_{m-2}\subset S_n(H_{m-2}, \pi_{m-2})$ and $\chi(S_n(H_{m-2}, \pi_{m-2}))\ge \chi(H_{m-2}) = m+1$ by Lemma~\ref{chrom nb = m+1}, and on the other hand, every vertex in $S_n(H_{m-2}, \pi_{m-2})\setminus H_{m-2}$ is of degree two in $S_n(H_{m-2}, \pi_{m-2})$, so every proper colouring of $H_{m-2}$ in $m+1\ge 3$ colours can be extended to a proper colouring of $\chi(S_n(H_{m-2}, \pi_{m-2}))$ in $m+1$ colours.
\end{proof}

Now, note that the random graph $PA_t(m, \delta)$ may be naturally seen as a digraph where every edge is oriented from its "younger" endvertex to the "older" one. Since the graph $PA_t(m, \delta)$ is increasing in $t$ with respect to inclusion, the probability to find a digraph $(H, \pi)$ as a subgraph of $PA_t(m, \delta)$ is also increasing in $t$ and therefore $\lim_{t\to +\infty} \mathbb P((H, \pi)\subset PA_t(m, \delta))$ is well defined. Our goal in this paper will be to prove that for every $m\ge 2$, the graph $H_{m-2}$ is realised as a subgraph of $PA_t(m, \delta)$ with probability, which tends to 1 as $t\to +\infty$. Our main tool will be the second moment method. However, by applying it directly to $H_{m-2}$ one can only obtain a lower bound on $\mathbb P((H_{m-2}, \pi_{m-2})\subset PA_t(m, \delta))$ that is a positive constant smaller than $1$. To improve this probability, we consider the sequence of subgraphs $(S_n(H_{m-2}, \pi_{m-2}))_{n\ge 1}$ instead.
\begin{theorem}\label{main th 1}
For all $m\ge 2$ and $\delta\in (-m,0)$, the sequence 
\begin{equation*}
(\lim_{t\to +\infty} \mathbb P(S_n(H_{m-2}, \pi_{m-2})\subset PA_t(m, \delta)))_{n\ge 1}  
\end{equation*}
converges to $1$.
\end{theorem}

\begin{proof}[Proof of Theorem~\ref{main th} assuming Theorem~\ref{main th 1}]
Theorem~\ref{main th} is a direct consequence of Observation~\ref{ob 2.4} and Theorem~\ref{main th 1}.
\end{proof}

To prove Theorem~\ref{main th 1} we analyse in detail the proof of Theorem 2.2 from \cite{GS} to compute the first and the second moment of the number of copies of the digraph $S_n(H_{m-2}, \pi_{m-2})$, equipped with a particular ordering, in $PA_t(m, \delta)$ for all $n\in \mathbb N$. Once having this information, we conclude by an application of the second moment method using Lemma~\ref{PZ ineq}.\par

To prepare the ground, we introduce some definitions from \cite{GS} and prove a couple of preliminary results. Fix $m\ge 2$ and $\delta\in (-m,0)$ and define $\tau = 3 + \dfrac{\delta}{m}\in (2,3)$. For a vertex $v$ in a weakly connected digraph $(H,\pi)$ define the quantity 
\begin{equation*}
    \beta(v)  = \beta_H(v) =  -\dfrac{\tau-2}{\tau-1} d^{(out)}_H(v) - \dfrac{1}{\tau - 1} d^{(in)}_H(v)
\end{equation*}
and let
\begin{align*}
    D(H, \pi) = & \max_{s = 0,1,\dots, |V(H)|} \left\{\sum_{i=s+1}^{|V(H)|} \left[1 -\dfrac{\tau-2}{\tau-1} d^{(out)}_H(\pi^{-1}(i)) - \dfrac{1}{\tau - 1} d^{(in)}_H(\pi^{-1}(i))\right]\right\}\\ = & \max_{s = 0,1,\dots, |V(H)|} \left\{\sum_{i=s+1}^{|V(H)|} (1+\beta(\pi^{-1}(i)))\right\}.
\end{align*}

For a digraph $(H, \pi)$ we say that $j$ is a \emph{maximiser} of $D(H, \pi)$ if
\begin{equation*}
    \sum_{i=j+1}^{|V(H)|} (1+\beta(\pi^{-1}(i))) = D(H, \pi).
\end{equation*}
Let the digraph $S_n(H, \pi)$ be equipped with the ordering $\sigma_{n, H}$ constructed above and assume that $H$ is a weakly connected digraph with at least two vertices.

\begin{lemma}\label{maximiser H_{m-2}}
For every $n > \tau - 1$, the unique maximiser of $D(S_n(H, \pi), \sigma_{n, H})$ is $s_H = |V(H)|$.
\end{lemma}
\begin{proof}
For every $n > \tau - 1$ and for every vertex $u\in V(H)$ we have that $d^{(in)}_{S_n(H, \pi)} (u)\ge n$, so
\begin{equation*}
    \beta_{S_n(H, \pi)}(u) = -\dfrac{\tau - 2}{\tau - 1}d^{(out)}_{S_n(H, \pi)}(u) - \dfrac{1}{\tau - 1} d^{(in)}_{S_n(H, \pi)}(u)\le -\dfrac{n}{\tau - 1} < -1,
\end{equation*}
and for every vertex $v\in V(S_n(H, \pi))\setminus V(H)$ we have that $d^{(out)}_{S_n(H,\pi)}(v)\ge 2$, so
\begin{equation*}
    \beta_{S_n(H, \pi)}(v) = -\dfrac{\tau - 2}{\tau - 1}d^{(out)}_{S_n(H, \pi)}(v) - \dfrac{1}{\tau - 1} d^{(in)}_{S_n(H, \pi)}(v) = - \dfrac{2(\tau - 2)}{\tau - 1} \in (-1, 0).
\end{equation*}
Since for every $u,v\in V(S_n(H, \pi))$ with $u\in V(H)$ and $v\in V(S_n(H, \pi))\setminus V(H)$ we have $\sigma_{n,H}(u) < \sigma_{n,H}(v)$ we conclude that the only maximiser of $D(S_n(H, \pi))$ is $|V(H)|$ independently of $n > \tau - 1$.
\end{proof}

From now on we denote for brevity $(H, \pi) = (H_{m-2}, \pi_{m-2}), s = |V(H)|\ge 7$ and $\beta_n(\cdot) = \beta_{S_n(H, \pi)}(\cdot)$. Note that $|V(S_n(H, \pi))| = s + \binom{s}{2} n$. For every $n\in \mathbb N$ we fix an ordering $\sigma_n = \sigma_{n,H}$ of $S_n(H, \pi)$ and denote $D_n = D(S_n(H, \pi), \sigma_n)$. 

By Lemma~\ref{maximiser H_{m-2}} we conclude that 
\begin{equation*}
    D_n = \sum_{i=s+1}^{s+\binom{s}{2}n} (1+\beta_n(\sigma_n^{-1}(i))) = \binom{s}{2} n \left(1-\dfrac{2(\tau-2)}{\tau-1}\right) = \dfrac{\binom{s}{2} n (3-\tau)}{\tau-1}.
\end{equation*}

\begin{figure}
\centering
\begin{tikzpicture}[line cap=round,line join=round,x=1cm,y=1cm]
\clip(-11.282695159292416,-3.5) rectangle (7.895934153996986,3.5);
\draw [->,line width=1.2pt] (-6,0) -- (-3,1.5);
\draw [->,line width=1.2pt] (-6,0) -- (-3,-1.5);
\draw [->,line width=1.2pt] (-6,0) -- (-4.5,0);
\draw [->,line width=1.2pt] (-3,0) -- (-1.5,1.5);
\draw [->,line width=1.2pt] (-3,0) -- (-1.5,-1.5);
\draw [->,line width=1.2pt] (-3,0) -- (-1.5,0);
\draw [->,line width=1.2pt] (0,3) -- (0,1.5);
\draw [->,line width=1.2pt] (0,-3) -- (0,-1.5);
\draw [->,line width=1.2pt] (-1,-1) -- (-0.5,-0.5);
\draw [->,line width=1.2pt] (-1,-1) -- (-2,-0.5);
\draw [->,line width=1.2pt] (-1,1) -- (-2,0.5);
\draw [->,line width=1.2pt] (-1,1) -- (-0.5,0.5);
\draw [->,line width=1.2pt] (1,1) -- (0.5,0.5);
\draw [->,line width=1.2pt] (1,1) -- (0.5,2);
\draw [->,line width=1.2pt] (1,-1) -- (0.5,-0.5);
\draw [->,line width=1.2pt] (1,-1) -- (0.5,-2);
\draw [line width=1.2pt] (-2,-0.5)-- (-3,0);
\draw [line width=1.2pt] (-0.5,-0.5)-- (0,0);
\draw [line width=1.2pt] (0,0)-- (-0.5,0.5);
\draw [line width=1.2pt] (-1.5,0)-- (0,0);
\draw [line width=1.2pt] (0,0)-- (0,1.5);
\draw [line width=1.2pt] (0.5,0.5)-- (0,0);
\draw [line width=1.2pt] (0,0)-- (0.5,-0.5);
\draw [line width=1.2pt] (0,0)-- (0,-1.5);
\draw [line width=1.2pt] (0.5,-2)-- (0,-3);
\draw [line width=1.2pt] (0,-3)-- (-1.5,-1.5);
\draw [line width=1.2pt] (-3,-1.5)-- (0,-3);
\draw [line width=1.2pt] (-3,1.5)-- (0,3);
\draw [line width=1.2pt] (-1.5,1.5)-- (0,3);
\draw [line width=1.2pt] (-4.5,0)-- (-3,0);
\draw [line width=1.2pt] (-2,0.5)-- (-3,0);
\draw [line width=1.2pt] (0.5,2)-- (0,3);
\begin{scriptsize}
\draw [fill=black] (-3,0) circle (2.5pt);
\draw [fill=black] (0,3) circle (2.5pt);
\draw [fill=black] (0,0) circle (2.5pt);
\draw [fill=black] (0,-3) circle (2.5pt);
\draw [fill=black] (-1,1) circle (2.5pt);
\draw [fill=black] (-1,-1) circle (2.5pt);
\draw [fill=black] (1,1) circle (2.5pt);
\draw [fill=black] (1,-1) circle (2.5pt);
\draw [fill=black] (-6,0) circle (2.5pt);

\draw [color=black] (-3.3,0.3) node {\Large $u_1$};
\draw [color=black] (0,3.3) node {\Large $u_2$};
\draw [color=black] (0.4,0) node {\Large $u_4$};
\draw [color=black] (0,-3.3) node {\Large $u_3$};
\draw [color=black] (-1,1.3) node {\Large $v_1$};
\draw [color=black] (-1,-1.3) node {\Large $v_2$};
\draw [color=black] (1.3,1) node {\Large $v_3$};
\draw [color=black] (1.3,-1) node {\Large $v_4$};
\draw [color=black] (-6.3,0) node {\Large $w$};
\end{scriptsize}
\end{tikzpicture}
\caption{In the figure there are two copies $S_1(H', \pi')$ and $S_1(H'', \pi'')$ of the digraph $S_1(H, \pi)$, where $(H, \pi)$ is the oriented triangle graph on the left in Figure~\ref{fig 1}, whose union forms a digraph with ordering $\hat{\sigma}_1$. These are induced by the sets of vertices $u_1, u_2, u_4, v_1, w, v_3$ and $u_1, u_3, u_4, v_2, w, v_4$ and may be equipped with the orderings $\sigma'_1: (u_4, u_2, u_1, w, v_1, v_3)\mapsto (1,2,3,4,5,6)$ and $\sigma''_1: (u_4, u_3, u_1, w, v_2, v_4)\mapsto (1,2,3,4,5,6)$. We may further define $\hat{\sigma}_1 : u_i\mapsto 4-i, w\mapsto 5, v_i\mapsto i+5$ for every $i\in [4]$. Indeed, we require that the ordering $\hat{\sigma}_1$ agrees with both $\sigma'_1$ and $\sigma''_1$ (assumption~\ref{cn 1}) and that $\hat{\sigma}_1, \sigma'_1, \sigma''_1$ satisfy $\hat{\sigma}_1(u_i) < \hat{\sigma}_1(w) < \hat{\sigma}_1(v_j)$ for all $i,j\in [4]$ (assumption~\ref{cn 2}).}
\label{fig 2}
\end{figure}
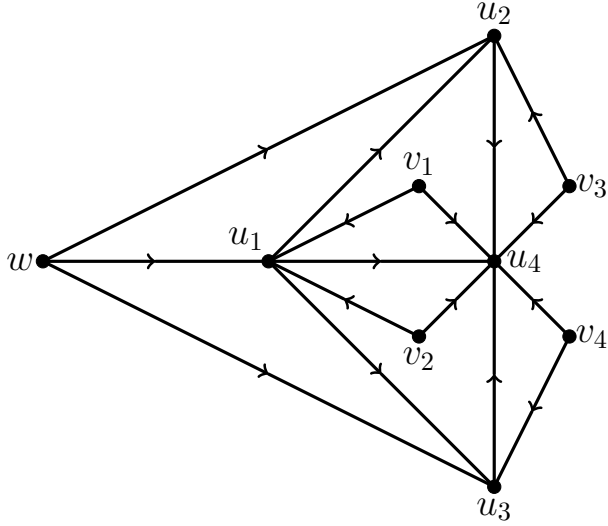

Let $(S_n(H', \pi'), \sigma'_n)$ and $(S_n(H'', \pi''), \sigma''_n)$ be two copies of the digraph with ordering $(S_n(H, \pi), \sigma_n)$, which may share common vertices and edges, where $n\ge 5$ is a fixed integer. We impose the restriction on $n$ to be sure that every vertex of degree at most four in $(S_n(H', \pi'), \sigma'_n)\cup (S_n(H'', \pi''), \sigma''_n)$ is outside $V(H'\cup H'')$. Moreover, assume that:
\begin{enumerate}
    \item\label{cn 1} $S_n(H', \pi')\cup S_n(H'', \pi'')$ is a digraph, which admits an ordering $\hat{\sigma}_n$ that agrees with both $\sigma'_n$ and $\sigma''_n$. In particular, the orientations of the common edges of $S_n(H', \pi')$ and $S_n(H'', \pi'')$, induced by $\sigma'_n$ and $\sigma''_n$ respectively, agree.
    \item\label{cn 2} for every vertex $v$ in $S_n(H', \pi')$ or in $S_n(H'', \pi'')$ of degree two in $S_n(H', \pi')\cup S_n(H'', \pi'')$ and for every vertex $w$ of degree two in $S_n(H', \pi')$ or in $S_n(H'', \pi'')$ and of degree at least three in $S_n(H', \pi')\cup S_n(H'', \pi'')$ we have $\hat{\sigma}_n(v) > \hat{\sigma}_n(w)$. See Figure~\ref{fig 2}.
\end{enumerate}

\noindent Denote $\hat{H}_n = S_n(H', \pi')\cup S_n(H'', \pi'')$ and $\hat{D}_n = D(S_n(H', \pi')\cup S_n(H'', \pi''), \hat{\sigma}_n)$.

\begin{lemma}
For every $n \ge 5$ we have that $\hat{D}_n\le 2D_n$ and equality holds if and only if $S_n(H', \pi')\cap S_n(H'', \pi'') = H'\cap H''$.
\end{lemma}
\begin{proof}
Remark that for the vertices $v$ of $H'\cup H''$ we have
\begin{equation*}
    \hat{\beta}_n(v) = \beta_{\hat{H}_n}(v) \le -\dfrac{n(s-1)}{\tau - 1} < -1.
\end{equation*}
On the other hand, the vertices of degree two in $S_n(H', \pi')$ and in $S_n(H'', \pi'')$ can be divided into three types.
\begin{enumerate}
    \item The vertices of type one are the ones, which appear in exactly one of the digraphs $S_n(H', \pi')$ and $S_n(H'', \pi'')$. For example, in Figure~\ref{fig 2} these are $v_1, v_2, v_3, v_4$. For any such vertex $v$ we have $\hat{\beta}_n(v) = -\dfrac{2(\tau - 2)}{\tau - 1}$. Denote the number of vertices of type one by $\ell_{1,n}$.
    \item The vertices of type two are the ones, which appear in each of the digraphs $S_n(H', \pi')$ and $S_n(H'', \pi'')$, but remain outside $H'\cup H''$. For example, in Figure~\ref{fig 2} the vertex $w$ is of type two. For any such vertex $v$ we have $\hat{\beta}_n(v) \le -\dfrac{3(\tau - 2)}{\tau - 1}$. Denote the number of vertices of type two by $\ell_{2,n}$
    \item The vertices of type three are the ones, which appear as degree two vertices in exactly one of the digraphs $S_n(H', \pi')$ and $S_n(H'', \pi'')$ and as vertices of degree at least three in the other digraph. Any such vertex $v$ participates in $H'\cup H''$ and therefore $\hat{\beta}_n(v) < -1$ by the above computation.
\end{enumerate}
Due to the assumptions \ref{cn 1} and \ref{cn 2} we have
\begin{align*}
    \hat{D}_n \le \hspace{0.25em}
    & \max \left\{\ell_{1,n}\left(1 - \dfrac{2(\tau - 2)}{\tau - 1}\right), \ell_{1,n}\left(1 - \dfrac{2(\tau - 2)}{\tau - 1}\right) + \ell_{2,n}\left(1 - \dfrac{3(\tau - 2)}{\tau - 1}\right)\right\}\\
    \le\hspace{0.25em}
    & (\ell_{1,n}+\ell_{2,n})\left(1 - \dfrac{2(\tau - 2)}{\tau - 1}\right) \le 2\binom{s}{2} n\left(1 - \dfrac{2(\tau - 2)}{\tau - 1}\right) = 2D_n.
\end{align*}
We have equality if and only if all vertices of degree two in $S_n(H', \pi')$ and in $S_n(H'', \pi'')$ are of type one or, otherwise said, if and only if the set of vertices of degree two in $\hat{H}_n$ is the union of the sets of vertices of degree two in $S_n(H', \pi')$ and in $S_n(H'', \pi'')$.
\end{proof}

\begin{corollary}\label{cor 2.8}
Each of the digraphs with ordering $(\hat{H}_n, \hat{\sigma}_n)$ obtained as union of two copies $S_n(H', \pi')$ and $S_n(H'', \pi'')$ of $S(H, \pi)$ with $\hat{D}_n = 2D_n$ may be constructed as follows:
\begin{enumerate}
    \item Construct the digraph $\hat{H}_0 = H'\cup H''$.
    \item Apply $S_n(\cdot)$ to the subgraph $H'$ of $\hat{H}_0$.
    \item Apply $S_n(\cdot)$ to the subgraph $H''$ of $\hat{H}_0$.\qed
\end{enumerate}
\end{corollary}

\begin{corollary}\label{cor 2.9}
The number of different digraphs $\hat{H}_n$ obtained as union of two copies $S_n(H', \pi')$ and $S_n(H'', \pi'')$ of $S_n(H, \pi)$ with $\hat{D}_n = 2D_n$ does not depend on $n$. Moreover, for every $n \ge 5$, each of these graphs has a unique maximiser of $\hat{D}_n$ given by $\hat{s} = s_{\hat{H}_n} = |V(H'\cup H'')|$.
\end{corollary}
\begin{proof}
The first statement is a direct consequence of Corollary~\ref{cor 2.8} - indeed, the number of digraphs $\hat{H}_n$ with $\hat{D}_n = 2D_n$ is equal to the number of digraphs, which may be constructed as union of two copies $(H', \pi')$ and $(H'', \pi'')$ of the digraph $(H, \pi)$ so that $\pi$ agrees with both orderings $\pi'$ and $\pi''$. For the second statement, by Corollary~\ref{cor 2.8} we have that
for every vertex $u\in H'\cup H''$ we have
\begin{equation*}
    \hat{\beta}_n(u) = -\dfrac{\tau - 2}{\tau - 1}d^{(out)}_{\hat{H}_n}(u) - \dfrac{1}{\tau - 1} d^{(in)}_{\hat{H}_n}(u)\le -\dfrac{n(s-1)}{\tau - 1} < -1
\end{equation*}
and for every vertex $v\in \hat{H}_n\setminus (H'\cup H'')$ we have
\begin{equation*}
    \hat{\beta}_n(v) = -\dfrac{\tau - 2}{\tau - 1}d^{(out)}_{\hat{H}_n}(v) - \dfrac{1}{\tau - 1} d^{(in)}_{\hat{H}_n}(v) = -\dfrac{2(\tau - 2)}{\tau - 1} \in (-1,0).
\end{equation*}
Since $\hat{\sigma}_n$ agrees with both $\sigma'$ and $\sigma''$ we know that for every vertex $u$ in $H'\cup H''$ and for every vertex $v$ in $\hat{H}_n\setminus (H'\cup H'')$ we have $\hat{\sigma}_n(u) < \hat{\sigma}_n(v)$. We conclude that $\hat{s} = |V(H'\cup H'')|$ is the only maximiser of $\hat{D}_n$. The corollary is proved.
\end{proof}

\section{Building upon the Garavaglia-Stegehuis result}\label{sec: GS result}

In this part, we mostly follow the notation of the paper \cite{GS}. Fix $m\in \mathbb N, m\ge 2$ and $\delta\in (-m,0)$. A digraph $(G, \sigma)$ is said to be \textit{attainable in $PA_t(m, \delta)$} if $(G, \sigma)$ could be realised as a subgraph of $PA_t(m, \delta)$ with positive probability. The following observation is straightforward.

\begin{observation}
The digraph $(G, \sigma)$ is attainable if $\max_{v\in V(G)} d^{(out)}_G(v) \le m$.\qed
\end{observation}

Now we state one of the main results of \cite{GS}. Fix an attainable digraph $(G, \sigma)$ and let $s_1, s_2, \dots, s_r\in [|V(G)|]$ be the maximisers of $D(G, \sigma)$. We denote by $N_t(G, \sigma)$ the number of copies of $(G, \sigma)$ in $PA_t(m,\delta)$ i.e. the number of occurrences of the digraph $G$ in $PA_t(m,\delta)$ so that the ordering $\sigma$ agrees with the canonical ordering of the vertices in $PA_t(m,\delta)$.

\begin{theorem}[\cite{GS}, Theorem 2.2]\label{GS Theorem}
For the attainable graph $(G, \sigma)$ with exactly $r$ maximisers of $D(G,\sigma)$ there are constants $C_1 = C_1(G), C_2 = C_2(G) > 0$ such that 
\begin{equation*}
    C_1\le \lim_{t\to +\infty} \dfrac{\mathbb E[N_t(G, \sigma)]}{t^{D(G, \sigma)}\log^{r-1}(t)}\le C_2.
\end{equation*}
\end{theorem}

Our approach requires a closer look in the proof of Theorem~\ref{GS Theorem}, which will provide more insight about the behaviour of the constants $C_1$ and $C_2$ as functions of the digraphs $(S_n(H, \pi), \sigma_n)_{n\ge 1}$. In particular, by Lemma~\ref{maximiser H_{m-2}} the unique maximiser $s = s_H = |V(H)|$ of $D_n$ remains fixed and independent of $n \ge 5$.\par

In the above setting the proof of Theorem~\ref{GS Theorem} from \cite{GS} establishes that there are absolute constants $C'_1, C'_2 > 0$, for which
\begin{align*}
    & C'_1 t^{D_n} \int_1^{\sqrt{t}} u_1^{\beta_n(\sigma_n^{-1}(1))} \int_{u_1}^{\sqrt{t}} u_2^{\beta_n(\sigma_n^{-1}(2))}\dots \int_{u_s}^{\sqrt{t}} u_s^{\beta_n(\sigma_n^{-1}(s))} du_s du_{s-1}\dots du_1\\
    & \cdot \int_{\frac{1}{\sqrt{t}}}^1 w_{s+1}^{\beta(\sigma_n^{-1}(s+1))}\int_{w_{s+1}}^1 w_2^{\beta(\sigma_n^{-1}(s+2))}\dots \int_{w_{s + \binom{s}{2} n -1}}^1 w_{s + \binom{s}{2} n}^{\beta(\sigma_n^{-1}(s + \binom{s}{2} n))} dw_{s + \binom{s}{2} n} dw_{s + \binom{s}{2} n -1}\dots dw_{s+1} \\
    \le\hspace{0.25em}
    & \mathbb E[N_t(S_n(H, \pi), \sigma_n)]\\
    \le\hspace{0.25em}
    & C'_2 t^{D_n} \int_1^{+\infty} u_1^{\beta_n(\sigma_n^{-1}(1))} \int_{u_1}^{+\infty} u_2^{\beta_n(\sigma_n^{-1}(2))}\dots \int_{u_s}^{+\infty} u_s^{\beta_n(\sigma_n^{-1}(s))} du_s du_{s-1}\dots du_1\\
    & \cdot \int_0^1 w_{s+1}^{\beta(\sigma_n^{-1}(s+1))}\int_{w_{s+1}}^1 w_2^{\beta(\sigma_n^{-1}(s+2))}\dots \int_{w_{s + \binom{s}{2} n -1}}^1 w_{s + \binom{s}{2} n}^{\beta(\sigma_n^{-1}(s + \binom{s}{2} n))} dw_{s + \binom{s}{2} n} dw_{s + \binom{s}{2} n -1}\dots dw_{s+1}.
\end{align*}

One has that for $t\to +\infty$, the limit of the expression
\begin{equation*}
    \int_1^{\sqrt{t}} u_1^{\beta_n(\sigma_n^{-1}(1))} \int_{u_1}^{\sqrt{t}} u_2^{\beta_n(\sigma_n^{-1}(2))}\dots \int_{u_s}^{\sqrt{t}} u_s^{\beta_n(\sigma_n^{-1}(s))} du_s du_{s-1}\dots du_1
\end{equation*}
is
\begin{equation*}
    \int_1^{+\infty} u_1^{\beta_n(\sigma_n^{-1}(1))} \int_{u_1}^{+\infty} u_2^{\beta_n(\sigma_n^{-1}(2))}\dots \int_{u_s}^{+\infty} u_s^{\beta_n(\sigma_n^{-1}(s))} du_s du_{s-1}\dots du_1 < \infty,
\end{equation*}
and the limit of the expression
\begin{equation*}
    \int_{\frac{1}{\sqrt{t}}}^1 w_{s+1}^{\beta(\sigma_n^{-1}(s+1))}\int_{w_{s+1}}^1 w_2^{\beta(\sigma_n^{-1}(s+2))}\dots \int_{w_{s+\binom{s}{2} n-1}}^1 w_{s+\binom{s}{2} n}^{\beta(\sigma_n^{-1}(s+\binom{s}{2} n))} dw_{s+\binom{s}{2} n} dw_{s+\binom{s}{2} n-1}\dots dw_{s+1}
\end{equation*}
is
\begin{equation*}
    \int_0^1 w_{s+1}^{\beta(\sigma_n^{-1}(s+1))}\int_{w_{s+1}}^1 w_2^{\beta(\sigma_n^{-1}(s+2))}\dots \int_{w_{s+\binom{s}{2} n-1}}^1 w_{s+\binom{s}{2} n}^{\beta(\sigma_n^{-1}(s+\binom{s}{2} n))} dw_{s+\binom{s}{2} n} dw_{s+\binom{s}{2} n-1}\dots dw_{s+1} < \infty.
\end{equation*}

Of course, these two limits depend on $n$. We have that
\begin{align*}
    & \int_1^{+\infty} u_1^{\beta_n(\sigma_n^{-1}(1))} \int_{u_1}^{+\infty} u_2^{\beta_n(\sigma_n^{-1}(2))}\dots \int_{u_s}^{+\infty} u_s^{\beta_n(\sigma_n^{-1}(s))} du_s du_{s-1}\dots du_1\\
    = \hspace{0.25em}
    & \int_1^{+\infty} u_1^{\beta_0(\sigma_0^{-1}(1)) - \frac{2n(\tau - 2)}{\tau - 1}} \int_{u_1}^{+\infty} u_2^{\beta_0(\sigma_0^{-1}(2))  - \frac{2n(\tau - 2)}{\tau - 1}}\dots \int_{u_s}^{+\infty} u_s^{\beta_0(\sigma_0^{-1}(s)) - \frac{2n(\tau - 2)}{\tau - 1}} du_s du_{s-1}\dots du_1 = \Theta\left(\dfrac{1}{n^s}\right).
\end{align*}

Moreover, since for every $i\ge s+1$ we have $\beta_n(\sigma^{-1}_n(i)) = -\dfrac{2(\tau-2)}{\tau-1}$,
\begin{align*}
    & \int_0^1 w_{s+1}^{\beta(\sigma_n^{-1}(s+1))}\int_{w_{s+1}}^1 w_{s+2}^{\beta(\sigma_n^{-1}(s+2))}\dots \int_{w_{s+\binom{s}{2} n-1}}^1 w_{s+\binom{s}{2} n}^{\beta(\sigma_n^{-1}(s+\binom{s}{2} n))} dw_{s+\binom{s}{2} n} dw_{s+\binom{s}{2} n-1}\dots dw_{s+1}\\
    =
    & \int_0^1 w_{s+1}^{-\frac{2(\tau-2)}{\tau-1}}\int_{w_{s+1}}^1 w_{s+2}^{-\frac{2(\tau-2)}{\tau-1}}\dots \int_{w_{s+\binom{s}{2} n-1}}^1 w_{s+\binom{s}{2} n}^{-\frac{2(\tau-2)}{\tau-1}} dw_{s+\binom{s}{2} n} dw_{s+\binom{s}{2} n-1}\dots dw_{s+1}
\end{align*}

Since $\tau\in (2,3)$, we know that the function $t\mapsto t^{-\frac{2(\tau-2)}{\tau-1}}$ is integrable on the interval $(0,1]$ and $\int_0^1 t^{-\frac{2(\tau-2)}{\tau-1}} dt = \dfrac{1}{1-\frac{2(\tau-2)}{\tau-1}} = \dfrac{\tau - 1}{3-\tau}$.

\begin{lemma}\label{lem 3.3}
\begin{equation*}
    \int_0^1 w_{s+1}^{-\frac{2(\tau-2)}{\tau-1}}\int_{w_{s+1}}^1 w_{s+2}^{-\frac{2(\tau-2)}{\tau-1}}\dots \int_{w_{s+\binom{s}{2} n-1}}^1 w_{s+\binom{s}{2} n}^{-\frac{2(\tau-2)}{\tau-1}} dw_{s+\binom{s}{2} n} dw_{s+\binom{s}{2} n-1}\dots dw_{s+1} = \left(\dfrac{\tau - 1}{3-\tau}\right)^{\binom{s}{2}n} \dfrac{1}{\left(\binom{s}{2}n\right)!}.
\end{equation*}
\end{lemma}
\begin{proof}
By summing the integrals 
\begin{equation*}
    \int_0^1 w_{\iota(s+1)}^{-\frac{\tau-2}{\tau-1}}\int_{w_{\iota(s+1)}}^1 w_{\iota(s+2)}^{-\frac{\tau-2}{\tau-1}}\dots \int_{w_{\iota(s+\binom{s}{2} n-1)}}^1 w_{\iota(s+\binom{s}{2} n)}^{-\frac{\tau-2}{\tau-1}} dw_{\iota(s+\binom{s}{2} n)} dw_{\iota(s+\binom{s}{2} n-1)}\dots dw_{\iota(s+1)}
\end{equation*}
for every permutation $\iota$ of $[s+\binom{s}{2} n]\setminus [s]$, we conclude by symmetry that the expression
\begin{equation*}
    \int_0^1 w_{s+1}^{-\frac{2(\tau-2)}{\tau-1}}\int_{w_{s+1}}^1 w_{s+2}^{-\frac{2(\tau-2)}{\tau-1}}\dots \int_{w_{s+\binom{s}{2} n-1}}^1 w_{s+\binom{s}{2} n}^{-\frac{2(\tau-2)}{\tau-1}} dw_{s+\binom{s}{2} n} dw_{s+\binom{s}{2} n-1}\dots dw_{s+1}
\end{equation*}
is a $\dfrac{1}{\left(\binom{s}{2}n\right)!}$-proportion of the integral
\begin{align*}
    & \int_0^1 w_{s+1}^{-\frac{2(\tau-2)}{\tau-1}}\int_0^1 w_{s+2}^{-\frac{2(\tau-2)}{\tau-1}}\dots \int_0^1 w_{s+\binom{s}{2} n}^{-\frac{2(\tau-2)}{\tau-1}} dw_{s+\binom{s}{2} n} dw_{s+\binom{s}{2} n-1}\dots dw_{s+1}\\ 
    =\hspace{0.3em} 
    & \left(\int_0^1 w_{s+1}^{-\frac{2(\tau-2)}{\tau-1}}\right)^{\binom{s}{2}n} =\hspace{0.3em} \left(\dfrac{\tau-1}{3-\tau}\right)^{\binom{s}{2}n}.
\end{align*}
The lemma is proved.
\end{proof}

\begin{corollary}\label{cor 3.4}
There are absolute constants $C''_1, C''_2 > 0$, for which, for every $n \ge 5$ and for every large enough $t$, we have
\begin{equation*}
    C''_1 t^{D_n} \left(\dfrac{\tau - 1}{3-\tau}\right)^{\binom{s}{2}n} \dfrac{1}{n^s \left(\binom{s}{2}n\right)!}\le \mathbb E[N_t(S_n(H, \pi), \sigma_n)]\le C''_2 t^{D_n} \left(\dfrac{\tau - 1}{3-\tau}\right)^{\binom{s}{2}n} \dfrac{1}{n^s \left(\binom{s}{2}n\right)!}.
\end{equation*}
\end{corollary}
\begin{proof}
This follows by Theorem~\ref{GS Theorem} and Lemma~\ref{lem 3.3}, using the fact that $D(S_n(H,\pi), \sigma_n)$ admits a unique maximiser by Lemma~\ref{maximiser H_{m-2}}.
\end{proof}

By analogous considerations we draw a similar conclusion for every digraph $(\hat{H}_n, \hat{\sigma}_n)$ obtained as union of two copies of $S_n(H, \pi)$, namely $S_n(H', \pi')$ and $S_n(H'', \pi'')$, and satisfying $\hat{D}_n = 2D_n$. Indeed, by Corollary~\ref{cor 2.9} for every $n\ge 5$ we have that $\hat{s} = s_{\hat{H}} = |V(H')\cup V(H'')|$ is the unique maximiser of $\hat{D}_n$.

\begin{corollary}\label{cor 3.5}
There are absolute constants $\hat{C}_1, \hat{C}_2 > 0$, for which, for every $n \ge 5$ and for every large enough $t$, we have
\begin{equation*}
    \hat{C}_1 t^{\hat{D}_n} \left(\dfrac{\tau - 1}{3-\tau}\right)^{2\binom{s}{2}n} \dfrac{1}{n^{\hat{s}} \left(2\binom{s}{2}n\right)!}\le \mathbb E[N_t(\hat{H}, \hat{\sigma}_n)]\le \hat{C}_2 t^{\hat{D}_n} \left(\dfrac{\tau - 1}{3-\tau}\right)^{2\binom{s}{2}n} \dfrac{1}{n^{\hat{s}} \left(2\binom{s}{2}n\right)!}.
\end{equation*}
\end{corollary}

\section{Proof of Theorem~\ref{main th 1}}\label{sec: proof}

By (\cite{GS}, Proof of Lemma 6.1) we know that, first, for any two copies $S_n(H', \pi')$ and $S_n(H'', \pi'')$ of $S_n(H, \pi)$ and for any two disjoint vertex sets $V', V''\subseteq (v_j)_{j\in [t]}$ we have
\begin{align*}
    & \mathbb P(\{S_n(H', \pi')\subseteq PA_t(m, \delta)[V']\} \cap \{S_n(H'', \pi'')\subseteq PA_t(m, \delta)[V'']\})\\
    =\hspace{0.25em}
    & \mathbb P(S_n(H', \pi')\subseteq PA_t(m, \delta)[V']) \mathbb P(\{S_n(H'', \pi'')\subseteq PA_t(m, \delta)[V'']).
\end{align*}

Therefore, computing the variance of $N_t(S_n(H, \pi), \sigma_n)$ yields
\begin{align}
    & \mathbb Var[N_t(S_n(H, \pi), \sigma_n)]\nonumber \\ 
    =\hspace{0.25em} 
    & \mathbb E[N_t(S_n(H, \pi), \sigma_n)^2] - \mathbb E[N_t(S_n(H, \pi), \sigma_n)]^2\nonumber\\
    \le\hspace{0.25em}
    &\label{line 1} \mathbb E[N_t(S_n(H, \pi), \sigma_n)] + \sum_{\substack{S_n(H', \pi'), S_n(H'', \pi'')\text{ copies of }S_n(H, \pi);\\H'\cap H'' \neq \varnothing}} \mathbb P((S_n(H', \pi')\cup S_n(H'', \pi''), \hat{\sigma}_n)\subseteq PA_t(m, \delta)).
\end{align}

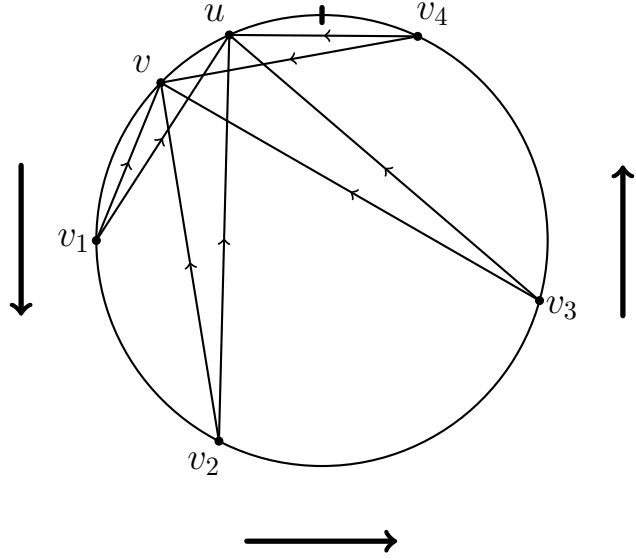
\begin{figure}
\centering
\begin{tikzpicture}[line cap=round,line join=round,x=1cm,y=1cm]
\clip(4.5,-4.25) rectangle (13.5,3.25);
\draw [line width=0.8pt] (9,0) circle (3cm);
\draw [->,line width=2pt] (5,1) -- (5,-1);
\draw [->,line width=2pt] (8,-4) -- (10,-4);
\draw [->,line width=2pt] (13,-1) -- (13,1);
\draw [line width=2pt] (9,2.9)-- (9,3.1);

\draw [->,line width=0.8pt] (6,0) -- (6.884018904503782, 1.3676868391526182);
\draw [line width=0.8pt] (6.884018904503782, 1.3676868391526182) -- (7.768037809007565,2.7353736783052365);

\draw [->,line width=0.8pt] (7.629807451145602,-2.668814789202107) -- (7.698922630076584, 0.033279444551565);
\draw [line width=0.8pt] (7.698922630076584, 0.033279444551565) -- (7.768037809007565,2.7353736783052365);

\draw [->,line width=0.8pt] (11.891588428755346,-0.7991973214971307) -- (9.829813118881456, 0.9680881784040529);
\draw [line width=0.8pt] (9.829813118881456, 0.9680881784040529) -- (7.768037809007565,2.7353736783052365);

\draw [->,line width=0.8pt] (10.272341231400265,2.716826787058545) -- (9.020189520203915,2.72610023268189075);
\draw [line width=0.8pt] (9.020189520203915,2.72610023268189075) -- (7.768037809007565,2.7353736783052365);

\draw [->,line width=0.8pt] (6,0) -- (6.42878573571858, 1.05);
\draw [line width=0.8pt] (6.42878573571858, 1.05) -- (9-2.14242852856285,2.1);

\draw [->,line width=0.8pt] (7.629807451145602,-2.668814789202107) -- (7.24368946129138, -0.284407);
\draw [line width=0.8pt] (7.24368946129138, -0.284407) -- (9-2.14242852856285,2.1);

\draw [->,line width=0.8pt] (11.891588428755346,-0.7991973214971307) -- (9.37457995009625, 0.650401);
\draw [line width=0.8pt] (9.37457995009625, 0.650401) -- (9-2.14242852856285,2.1);

\draw [->,line width=0.8pt] (10.272341231400265,2.716826787058545) -- (8.56495635141871, 2.40841);
\draw [line width=0.8pt] (8.56495635141871, 2.40841) -- (9-2.14242852856285,2.1);

\begin{scriptsize}
\draw [fill=black] (7.768037809007565,2.7353736783052365) circle (1.5pt);
\draw [fill=black] (9-2.14242852856285,2.1) circle (1.5pt);
\draw [fill=black] (6,0) circle (1.5pt);
\draw [fill=black] (7.629807451145602,-2.668814789202107) circle (1.5pt);
\draw [fill=black] (11.891588428755346,-0.7991973214971307) circle (1.5pt);
\draw [fill=black] (10.272341231400265,2.716826787058545) circle (1.5pt);

(9-2.14242852856285,2.1)

\draw[color=black] (7.768037809007565-0.2,2.7353736783052365+0.3) node {\Large $u$};
\draw[color=black] (9-2.14242852856285-0.25,2.1+0.25) node {\Large $v$};
\draw[color=black] (6-0.3,0) node {\Large $v_1$};
\draw[color=black] (7.629807451145602-0.2,-2.668814789202107-0.3) node {\Large $v_2$};
\draw[color=black] (11.891588428755346+0.3,-0.7991973214971307-0.1) node {\Large $v_3$};
\draw[color=black] (10.272341231400265+0.2,2.716826787058545+0.3) node {\Large $v_4$};
\end{scriptsize}
\end{tikzpicture}
\caption{In the figure $n=2$ and $u,v\in H'\cap H''$ have four common neighbours $v_1, v_2, v_3, v_4$ of degree two in $\hat{H}_2$. The arrows outside the circle show the order, in which the vertices in the PA graph are added, and the short thick segment on the top of the circle marks the beginning and the end. From $v_1, v_2, v_3, v_4$, any subset of two vertices could have come from $S_2(H', \pi')$ (and the other two vertices respectively from $S_2(H'', \pi'')$). The order of the vertices in every group in the ordering $\hat{\sigma}_2$ fixed in advance as $\hat{\sigma}_2$ agrees with both $\sigma'_2$ and $\sigma''_2$.}
\label{fig 3}
\end{figure}

Now, since the number of unions of two copies $(H', \pi')$ and $(H'',\pi'')$ of $(H, \pi)$ is finite, by Corollary~\ref{cor 2.9} the number of graphs $\hat{H}_n$ with $\hat{D}_n = 2D_n$ is finite and does not depend on $n \ge 5$. However, for fixed $(H', \pi')$ and $(H'', \pi'')$ the graph $(\hat{H}_n, \hat{\sigma}_n)$ may be obtained as union of two copies $(S_n(H', \pi'), \sigma'_n)$ and $(S_n(H'', \pi''), \sigma''_n)$ of the digraph $(S_n(H, \pi), \sigma_n)$ in different ways. Indeed, for every pair of vertices $u,v$ in $H'\cap H''$, the digraph $\hat{H}_n$ contains $2n$ vertices of degree two incident to oriented edges to $u$ and $v$ and any subset of these $2n$ vertices could come from $S_n(H', \pi')$ a priori. Of course, this observation may only lead to an upper bound due to the condition that $\hat{\sigma}_n$ must agree with both $\sigma'_n$ and $\sigma''_n$. See Figure~\ref{fig 3}.

Therefore, one has that the number of different pairs of copies $(S_n(H', \pi'), \sigma'_n)$ and $(S_n(H'', \pi''), \sigma''_n)$ of the digraph $(S_n(H, \pi), \sigma_n)$, which form the same graph $\hat{H}_n$, is at most
\begin{equation*}
   \dbinom{2n}{n}^{\binom{|V(H'\cap H'')|}{2}} = \dbinom{2n}{n}^{\binom{2s-\hat{s}}{2}}. 
\end{equation*}
By Corollary~\ref{cor 2.9} this gives an upper bound for (\ref{line 1}) that yields
\begin{align}
    & \mathbb E[N_t(S_n(H, \pi), \sigma_n)]\nonumber \\ +\hspace{0.25em}
    & (1+o_t(1)) \sum_{\substack{\hat{H}_n: \text{ there exist}\\ S_n(H', \sigma'), S_n(H'', \sigma'')\text{ copies of }S_n(H, \pi);\\ S_n(H', \sigma')\text{ and } S_n(H'', \sigma'') \text{ have no common degree two vertices};\\ \hat{H}_n = S_n(H', \pi')\cup S_n(H'', \pi'')\text{ and } H'\cap H'' \neq \varnothing}} \dbinom{2n}{n}^{\binom{2s-\hat{s}}{2}} \mathbb E((\hat{H}_n, \hat{\sigma}_n)\subseteq PA_t(m, \delta))\nonumber \\
    \le\hspace{0.25em}
    & \mathbb E[N_t(S_n(H, \pi), \sigma_n)]\nonumber \\ +\hspace{0.25em} 
    &\label{line 2} (1+o_t(1)) \sum_{\substack{\hat{H}_n: \text{ there exist}\\ S_n(H', \sigma'), S_n(H'', \sigma'')\text{ copies of }S_n(H, \pi);\\ S_n(H', \sigma')\text{ and } S_n(H'', \sigma'') \text{ have no common degree two vertices};\\ \hat{H}_n = S_n(H', \pi')\cup S_n(H'', \pi'')\text{ and } H'\cap H'' \neq \varnothing}} \dbinom{2n}{n}^{\binom{2s-\hat{s}}{2}} \hat{C}_2 t^{\hat{D}_n} \left(\dfrac{\tau - 1}{3-\tau}\right)^{2\binom{s}{2}n} \dfrac{1}{n^{\hat{s}} \left(2\binom{s}{2}n\right)!}.
\end{align}

Moreover, since $H'\cap H''\neq \varnothing$ we have that $s\le \hat{s} = |V(H'\cup H'')| < 2|V(H)| = 2s$. By direct application of Stirling's formula we conclude that
\begin{align}
    & \dbinom{2n}{n}^{\binom{2s-\hat{s}}{2}}\left(\dfrac{\tau - 1}{3-\tau}\right)^{2\binom{s}{2}n} \dfrac{1}{n^{\hat{s}} \left(2\binom{s}{2}n\right)!} = \dbinom{2n}{n}^{\binom{2s-\hat{s}}{2}} n^{2s-\hat{s}} \dfrac{(\binom{s}{2} n)!^2}{(2\binom{s}{2} n)!} \left(\left(\dfrac{\tau - 1}{3-\tau}\right)^{\binom{s}{2}n} \dfrac{1}{n^s \left(\binom{s}{2}n\right)!}\right)^2\nonumber\\ 
    =\hspace{0.25em} 
    & \Theta_n\left(2^{\left(2\binom{2s-\hat{s}}{2}-2\binom{s}{2}\right)n} n^{\frac{2(2s-\hat{s})+1-\binom{2s-\hat{s}}{2}}{2}}\right) \left(\left(\dfrac{\tau - 1}{3-\tau}\right)^{\binom{s}{2}n} \dfrac{1}{n^s \left(\binom{s}{2}n\right)!}\right)^2\nonumber\\
    = \hspace{0.25em}
    &\label{line 3} o_n(1) \left(\left(\dfrac{\tau - 1}{3-\tau}\right)^{\binom{s}{2}n} \dfrac{1}{n^s \left(\binom{s}{2}n\right)!}\right)^2.
\end{align}

The last equality comes from the fact that if $\hat{s} > s$, $2^{\left(2\binom{2s-\hat{s}}{2}-2\binom{s}{2}\right)n}$ decreases exponentially in $n$, and if $\hat{s} = s$ (in which case $H'$ and $H''$ coincide) then
\begin{equation*}
2(2s-\hat{s})+1-\binom{2s-\hat{s}}{2} = 2s+1-\dfrac{s(s-1)}{2} = 1 + \dfrac{s(5-s)}{2} < 0
\end{equation*}
for every $s\ge 7$. Here we use the fact that for every $i\ge 0$ we have $|V(H_i)|\ge |V(H_0)| = 7$. 

Thus, since by Corollary~\ref{cor 2.9} the number of terms in sum in (\ref{line 2}) depends only on $s$ and not on $t$ and $n$, by Corollary \ref{cor 3.4}, (\ref{line 2}) and (\ref{line 3}) we conclude that the limit
\begin{align*}
    \lim_{t\to +\infty} \dfrac{\mathbb E[N_t(S_n(H, \pi), \sigma_n)]^2}{\mathbb E[N_t(S_n(H, \pi), \sigma_n)^2]} = \lim_{t\to +\infty} \dfrac{\mathbb E[N_t(S_n(H, \pi), \sigma_n)]^2}{\mathbb E[N_t(S_n(H, \pi), \sigma_n)]^2 + \mathbb Var[N_t(S_n(H, \pi), \sigma_n)]}
\end{align*}
exists and tends to $1$ as $n\to +\infty$. This is sufficient to conclude by Lemma~\ref{PZ ineq} that
\begin{equation*}
    \lim_{t\to +\infty} \mathbb P(S_n(H, \pi)\subset PA_t(m, \delta)) = \lim_{t\to +\infty} \mathbb P(N_t(S_n(H, \pi)) > 0)
\end{equation*}
converges to $1$ as $n\to +\infty$. Theorem~\ref{main th 1} is proved.

\section{Conclusion and open problems}\label{sec: conclusion}
In this paper we prove that the chromatic number of the random graph $PA_t(m, \delta)$ converges almost surely to $m+1$ when $t\to +\infty$ for every $m\ge 1$ and $\delta\in (-m,0)$. To do this, we construct a family of digraphs $(S_n(H, \pi))_{n\ge 1}$ and prove that some graph in this family is present as a subgraph of $PA_t(m, \delta)$ asymptotically almost surely. It is a natural question to study the speed of convergence to this almost sure limit. One may also wonder what happens if $\tau \ge 3$. We conjecture that the same phenomenon takes place.
\begin{conjecture}
For every $m\ge 1$ and every $\delta > -m$, the chromatic number of $(PA_t(m, \delta))_{t\ge 1}$ converges almost surely to $m+1$.
\end{conjecture}

\section{Acknowledgements} 
The author would like to thank Dieter Mitsche and Violeta Naydenova for several corrections. Thanks are also due to the two anonymous referees for a number of useful comments and suggestions.

\bibliographystyle{plain}
\bibliography{References}

\end{document}